\documentclass[a4paper,11pt]{article}
\usepackage{url}
\usepackage{amssymb}
\usepackage{amsthm}
\usepackage{amsmath}
\usepackage{anysize}
\usepackage{graphicx}
\usepackage{subfigure}
\usepackage{psfrag}

\newcommand\R{\mathbb{R}}
\newcommand\Rd[1]{(\mathbb{R}^{#1})^*}
\newcommand\aff{\mathsf{aff}\,}
\newcommand\dual{\Delta}
\newcommand\fdual{\diamond}
\newcommand\K{\mathsf{K}}
\newcommand\Ll{\mathsf{L}}
\newcommand\cc[1]{{#1}^{\mathsf{c}}}
\newcommand\bd{\mathcal{B}}
\renewcommand\S{\mathbb{S}}
\renewcommand\dim{\mathsf{dim}\,}
\newcommand\Zt{{\mathbb{Z}_2}}
\newcommand\Ztmap{\rightarrow_\Zt}
\newcommand\djn[1]{{#1}^{*2}_\Delta}
\newcommand\ind{\mathsf{ind}_{\Zt}\,}
\newcommand\KG{\mathsf{KG}}
\newcommand\nf{\mathcal{F}}
\renewcommand\ker{\mathsf{ker}\,}
\newcommand\im{\mathsf{im}\,}
\renewcommand\int{\mathsf{int}\,}
\newcommand\nc{\mathcal{N}}
\newcommand\V{\mathsf{vert}\,}
\newcommand\conv{\mathsf{conv}\,}
\newcommand\cone{\mathsf{cone}\,}
\newcommand\A{\mathcal{A}}
\newcommand\C{\mathcal{C}}
\newcommand\eps{\varepsilon}

\newtheorem{thm}{Theorem}[section]
\newtheorem{cor}[thm]{Corollary}
\newtheorem{lem}[thm]{Lemma}
\newtheorem{prop}[thm]{Proposition}

\theoremstyle{definition}
\newtheorem{dfn}[thm]{Definition}

\newtheorem{obs}{Observation}

\begin{document}
\title{Topological obstructions for vertex numbers of Minkowski sums}

\author{Raman Sanyal\footnote{%
This work was supported by the Deutsche Forschungsgemeinschaft via a Leibniz
Grant and the research training group ``Methods for Discrete Structures'' (GRK
1408).}\\\small 
Institut f\"ur Mathematik, MA 6-2 \\  \small TU Berlin, 10623 Germany\\
\small sanyal@math.tu-berlin.de}

\date{March 22, 2007}

\maketitle

\begin{abstract}\noindent
We show that for polytopes $P_1, P_2, \dots, P_r \subset \R^d$, each having
$n_i \ge d+1$ vertices, the Minkowski sum $P_1 + P_2 + \cdots + P_r$ cannot
achieve the maximum of $\prod_i n_i$ vertices if $r \ge d$.  This complements
a recent result of Fukuda \& Weibel (2006), who show that this is possible for
up to $d-1$ summands. The result is obtained by combining methods from
discrete geometry (Gale transforms) and topological combinatorics (van
Kampen--type obstructions) 
\end{abstract}

\section{Introduction}\label{sec:intro}

For two polytopes $P, Q \subset \R^d$ their \emph{Minkowski sum} is the
(innocent--looking) convex polytope
\[
    P + Q = \{ p + q : p \in P,\, q \in Q \} \subset \R^d.
\]
Minkowski sums have starred in applied areas, such as robot motion
planning \cite{latombe} and computer aided design
\cite{DBLP:journals/cad/ElberK99}, as well as in fields of pure mathematics,
among them commutative algebra and tropical geometry \cite{stu02}.  In
applications it is essential to understand the facial structure of $P + Q$.
But, even with quite detailed knowledge of $P$ and $Q$, it is in general
difficult to determine the combinatorics of $P+Q$.  Even for {\it special
cases}, the knowledge of complete face lattices is meager. The best understood
Minkowski sums are zonotopes \cite{zie95} and sums of perfectly centered
polytopes with their polar duals \cite{fuk06}.

So it is natural (and vital) to investigate the combinatorial structure of
Minkowski sums. 
From the standpoint of combinatorial geometry, a less ambitious goal one can
settle for is the question of $f$-vector shapes. This includes different kinds
of upper and lower bounds for the $f$-vector entries with respect to the
corresponding entries of the summands. Starting with the first entry of an
$f$-vector, the number of vertices, Fukuda \& Weibel \cite{fuk06} studied the
maximal number of vertices of a Minkowski sum. Their starting point was the
following upper bound on the number of vertices.

\begin{prop}[Trivial Upper Bound; cf.\ \cite{fuk06}]\label{prop:triv_bound}
    Let $P_1,\dots,P_r \subset \R^d$ be polytopes. Then
    \[
        f_0(P_1 + \cdots + P_r) \ \ \le \ \ \prod_{i=1}^r f_0(P_i).
    \]
\end{prop}

Fukuda and Weibel gave a construction that showed that the trivial upper bound
can be attained independent of the dimension but with a restricted number of
summands.

\begin{thm}[Fukuda \& Weibel \cite{fuk06}]
    For every $r < d$ there are $d$-polytopes $P_1,\dots,P_r \subset \R^d$,
    each with arbitrarily large number of vertices, such that the Minkowski
    sum $P_1 + \cdots + P_r$ attains the trivial upper bound on the number of
    vertices.
\end{thm}

This result was our point of departure.  In this paper, we set out to prove
the following result that asserts that the restriction to the number of
summands is best possible.

\begin{thm}\label{thm:main}
    Let $r \ge d$ and let $P_1,\dots,P_r \subset \R^d$ be polytopes with
    $f_0(P_i) \ge d+1$ vertices for all $i = 1,\dots, r$.  Then
    \[
        f_0(P_1 + \cdots + P_r) \ \le\ \left( 1 -
        \frac{1}{(d+1)^r}\right)\prod_{i=1}^r f_0(P_i).
    \]
\end{thm}

However, based on a preliminary version of this paper, C.\ Weibel (personal
communication) noted that Theorem \ref{thm:main} is not optimal with respect
to the number of vertices of the summands.

The paper is organized in the following manner: In Section
\ref{sec:problem_reductions} we gather some observations that will reduce the
statement to one special case (per dimension) whose validity has to be
checked. In particular, we give a reformulation of the problem that casts it
into a stronger question concerning projections of polytopes.  The very
structure of the problem at hand allows for the utilization of the tools
devised in R\"{o}rig, Sanyal \& Ziegler \cite{rsz07}.  In Section
\ref{sec:proj_properties} we give a review to the terminology and techniques.

The punchline will be that a realization of a polytope with certain properties
under projection gives rise to, first, a polytope associated to the projection
and, secondly, to a simplicial complex that is embedded in the boundary of
this polytope. To prove the non-existence of a certain realization it will
suffice to show that this complex is not embeddable into a sphere of the
prescribed dimension. To show the reader that we did not deal one difficult
problem for another one, we give, in Section \ref{sec:embeddability}, a short
account of Matou\v{s}ek's book \cite{MatousekBZ:BU}, in which he presents
means for dealing with embeddability questions in a combinatorial fashion. In
Section \ref{sec:compl_complex}, we finally put together the pieces gathered
to prove a stronger version of Theorem \ref{thm:main} and we close with some
remarks in Section \ref{sec:remarks}. 

{\bf Acknowledgements}. 
We would like to express our gratitude to Christian Haase, Thilo R\"orig, and
G\"unter M. Ziegler for stimulating discussions and valuable insights. We
also would like to thank Bernd Sturmfels for helpful comments on the exposition.

\section{The problem, some reductions \& a reformulation}
\label{sec:problem_reductions}

Forming Minkowski sums is not a purely combinatorial construction, i.e.\ in
contrast to basic polytope constructions such as products, direct sums, joins,
etc.\ the resulting face lattice is not determined by the face lattices of the
polytopes involved.  For a sum $P +Q$ of two polytopes $P$ and $Q$ it is easy
to see that if $F \subseteq P + Q$ is a proper face, then $F$ is of the form
$F = G + H$ with $G \subseteq P$ and $H \subseteq Q$ being faces.  This, in
particular, sheds new light on the ``Trivial Upper Bound'' in the last
section: It states that the set of vertices of a Minkowski sum is a subset of
the pairwise sums of vertices of the polytopes involved.

As a guiding example let us consider the first non-trivial case: Are there two
triangles $P$ and $Q$ in the plane whose sum is a $9$-gon? 
An \emph{ad-hoc} argument  for this case, that uses notation and terminology
presented in \cite{zie95}, is the following: Clearly, the polytope $P+Q$ is a
$9$-gon if its normal fan $\nc(P+Q)$ has nine extremal rays. The normal fan
$\nc(P+Q)$ equals $\nc(P) \wedge \nc(Q)$, the common refinement of the fans
$\nc(P)$ and $\nc(Q)$. Thus, the cones in $\nc(P+Q)$ are pairwise
intersections of cones of $\nc(P)$ and $\nc(Q)$. It follows that the extremal
rays, i.e.\ the $1$-dimensional cones, of $\nc(P+Q)$ are just the extremal
rays of $P$ and of $Q$. If $P$ and $Q$ are triangles, then each one has only
three extremal rays. Therefore, $\nc(P+Q)$ has at most six extremal rays and
falls short of being a $9$-gon. The same reasoning yields the following
result.

\begin{prop} 
    Let $P$ and $Q$ be two polygons in the plane. Then
    \[
        f_0(P+Q) \le f_0(P) + f_0(Q).
    \]
\end{prop}

However, this elementary geometrical reasoning fails in higher dimensions and
we will employ topological machinery for the general case. But for now let us
give some observations that will simplify the general case.

The first observation concerns the dimensions of the polytopes involved in the
sum. 

\begin{obs}[Dimension of summands]
    Let $P_1, P_2, \dots,P_r \subset \R^d$ be polytopes each having at least
    $d+1$ vertices. Then there are \emph{full-dimensional} polytopes
    $P^\prime_1,P^\prime_2, \dots,P^\prime_r$ with $f_0(P^\prime_i) =
    f_0(P_i)$ and $f_0(P^\prime_1 + P^\prime_2 +  \cdots + P^\prime_r) \ge
    f_0(P_1 + P_2 +  \cdots + P_r)$
\end{obs}

Clearly, if one of the summands, say $P_1$, is not full-dimensional, then the
number of vertices prevents $P_1$ from being a lower dimensional simplex.
Choosing a vertex $v$ of $P_1$ that is not a cone point and pulling $v$ in a
direction perpendicular to its affine hull yields a polytope $P^\prime_1$ with
$f_0(P_1) = f_0(P^\prime_1)$ with $\dim P^\prime_1 = \dim P + 1$.  Exchanging
$P_1$ for $P^\prime_1$ possibly increases the number of vertices of the
Minkowski sum.

\begin{obs}[Number of summands] 
    Let $P_1, P_2,\dots,P_r \subset \R^d$ be $d$-polytopes such that $P_1 +
    P_2 + \cdots + P_r$ attains the trivial upper bound, then so does every
    subsum $P_{i_1} + P_{i_2} +  \cdots + P_{i_k}$ with $\{ i_1, \dots i_k \}
    \subseteq [r]$.
\end{obs}

Thus we can restrict to the situation of sums with $d$ summands.  The next
observation turns out to be even more valuable. It states that we can even
assume that every summand is a simplex. 

\begin{obs}[Combinatorial type of summands]
    Let $P_1, P_2,\dots,P_r \subset \R^d$ be $d$-polytopes such that $P_1 +
    P_2 + \cdots + P_r$ attains the trivial upper bound. For every $i \in [r]$
    let $P^\prime_i \subset P_i$ be a vertex induced, full-dimensional
    subpolytope, then $P^\prime_1 + P^\prime_2 + \cdots + P^\prime_r$ attains
    the trivial upper bound.
\end{obs}

The last observation casts the problem into the realm of polytope projections. 

\begin{obs} 
    The Minkowski sum $P + Q$ is the projection of the product $P \times Q$
    under the map $\pi: \R^d \times \R^d \rightarrow \R^d$ with $ \pi(x,y) = x
    + y$.
\end{obs}

We will derive  Theorem \ref{thm:main} from the following stronger statement.

\begin{thm} \label{thm:main_stronger}
    Let $P$ be a polytope combinatorially equivalent to a $d$-fold product of
    $d$-simplices and let $\pi : P \rightarrow \R^d$ be a linear projection.
    Then
    \[
        f_0(\pi P) \le f_0(P) - 1 = (d+1)^d - 1.
    \]
\end{thm}

Before we give a proof of Theorem \ref{thm:main} let us remark on a few things
concerning the previous theorem.

Special emphasis should be put on the phrase ``combinatorially equivalent to a
product.'' The \emph{standard} product  $P \times Q$ of two polytopes
$P\subset \R^d$ and $Q\subset \R^e$ is obtained by taking the Cartesian
product of $P$ and $Q$, that is taking the convex hull of $\V P \times \V Q
\subset \R^{d+e}$.  One feature of this construction is that if $P^\prime
\subset P$ and $Q^\prime \subset Q$ are vertex induced subpolytopes, then
$P\times Q$ contains $P^\prime \times Q^\prime$ again as a vertex induced
subpolytope. This no longer holds for \emph{combinatorial products}: The
\emph{Goldfarb cube} $G_4$ is a $4$-polytope combinatorially equivalent to a
$4$-cube with the property that a projection to $2$-space retains all $16$
vertices (cf.\ \cite{Goldfarb2,am98}). The $4$-cube itself is combinatorially
equivalent to a product of two quadrilaterals which, in the standard product,
contains a vertex induced product of two triangles.  The subpolytope on the
corresponding vertices of $G_4$ is not a combinatorial product of two
triangles; indeed, this would contradict Theorem \ref{thm:main_stronger} for
$d=2$.

The bound of Theorem \ref{thm:main_stronger} seems to be tight: In Section
\ref{sec:remarks} we give a realization of a product of two triangles such
that a projection to $2$-space has $8$ vertices. Therefore, the bound given in
Theorem \ref{thm:main} is \emph{not} tight for $d = 2$: The sum of two
triangles in the plane has at most $6$ vertices. 

Theorem \ref{thm:main_stronger} touches upon properties of the realization
space of products of simplices. While, in general, realization spaces are
rather delicate objects, the statement at hand is on par with the fact that
positively spanning vector configurations with prescribed sign patterns of
linear dependencies do not exist. Methods for treating such problems were
developed in \cite{rsz07}; the next two sections give a made-to-measure
introduction.

\begin{proof}[Proof of Theorem \ref{thm:main}]
    We can assume that all polytopes $P_i$ have their vertices $V_i$ in
    general position. For every choice $V^\prime_i \subseteq V_i$ of $d+1$
    vertices from each $P_i$ the polytopes $P^\prime_i = \conv V^\prime_i$ are
    $d$-simplices. The Minkowski sum $P^\prime_1+ P^\prime_2+\cdots
    +P^\prime_r$ is a projection of a product of $r \ge d$ simplices and, by
    Theorem \ref{thm:main_stronger}, does not attain the trivial upper bound.
    There are exactly $\prod_{i=1}^r \tbinom{f_0(P_i)}{d+1}$ choices for the
    $P^\prime_i$. On the other hand, every sum of vertices
    $v_1+v_2+\cdots+v_r$ occurs in only $\prod_{i=1}^r\tbinom{f_0(P_i)-1}{d}$
    different subsums.  Thus, by the pigeonhole principle, there are at least 
    \[
        \prod_{i=1}^r \frac{\tbinom{f_0(P_i)}{d+1}}{\tbinom{f_0(P_i)-1}{d}} =
        \prod_{i=1}^r \frac{f_0(P_i)}{d+1}
    \]
    sums of vertices that fail to be a vertex in at least one subsum.
\end{proof}

\section{Geometric and Combinatorial Properties of Projections}
\label{sec:proj_properties}

In polytope projections faces can collapse or get mapped to the interior.
Therefore, it is difficult to predict the (combinatorial) outcome of a
projection. There is, however, a class of faces that behave nicely under
projection and whose properties we will exploit in the following.

\begin{dfn}[Strictly preserved faces, cf. \cite{zie04}]
    Let $P$ be a polytope, $F \subseteq P$ a proper face and $\pi: P
    \rightarrow \pi(P)$ a linear projection of polytopes. The face $F$ is
    \emph{strictly preserved} under $\pi$ if 
    \begin{enumerate}
        \item[  i)] $H = \pi(F)$ is a face of $\pi(P)$,
        \item[ ii)] $F$ and $H$ are combinatorially isomorphic, and
        \item[iii)] $\pi^{-1}(H)$ is equal to $F$.
    \end{enumerate}
\end{dfn}

\begin{figure}[ht]
    \centering
    \psfrag{p}[bl]{$\scriptstyle \pi$}
    \psfrag{q}[cc]{$\scriptstyle q^*$}

    \psfrag{v1}[bl]{$\scriptstyle v_1$}
    \psfrag{v2}[tl]{$\scriptstyle v_2$}
    \psfrag{v3}[cr]{$\scriptstyle v_3$}
    \psfrag{v4}[br]{$\scriptstyle v_4$}
    \psfrag{v5}[bl]{$\scriptstyle v_5$}

    \psfrag{pv1}[cc]{$\scriptstyle \pi(v_1)$}
    \psfrag{pv2}[cc]{$\scriptstyle \pi(v_2)=\pi(v_5)$}
    \psfrag{pv3}[cc]{$\scriptstyle \pi(v_3)=\pi(v_4)$}

    \psfrag{dv1}[cl]{$\scriptstyle v_1^\fdual$}
    \psfrag{dv2}[cb]{$\scriptstyle v_2^\fdual$}
    \psfrag{dv3}[cb]{$\scriptstyle v_3^\fdual$}

    \psfrag{qv1}[tr]{$\scriptstyle q^*(v_1^\fdual)$}
    \psfrag{qv2}[cr]{$\scriptstyle q^*(v_2^\fdual)$}
    \psfrag{qv3}[cr]{$\scriptstyle q^*(v_3^\fdual)$}
    \includegraphics[height=6cm]{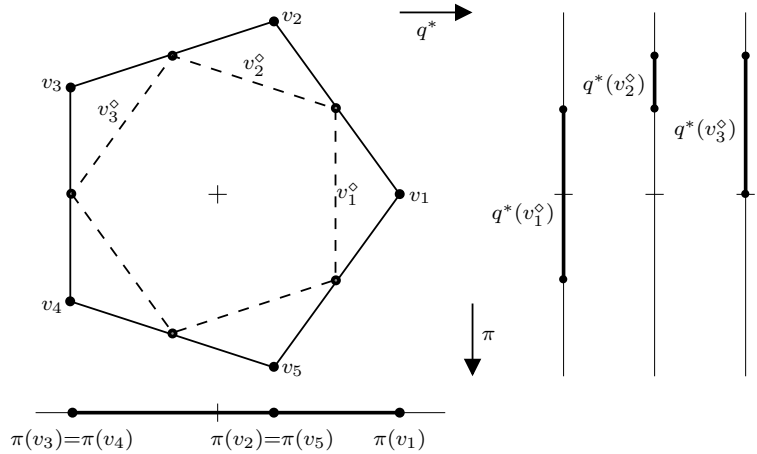}
    \caption{\label{fig:faces}Projection of pentagon to the real line: $v_1$
    is strictly preserved, $v_2$ is not preserved, and $v_3$ is preserved but
    not strictly. In accordance with the Projection Lemma: $q^*(v_1^\fdual)$
    contains $0$ in interior, $q^*(v_2^\fdual)$ does not, and $q^*(v_3^\fdual)$
    has $0$ in its boundary.}
\end{figure}

The first two conditions should trigger an agreeing nod  since they model the
intuition behind ``preserved faces.'' The third condition is a little less
clear. In order to talk about \emph{distinct} faces of a projection we have to
rule out that two preserved faces come to lie on top of each other and this
situation is dealt with in condition iii). Figure \ref{fig:faces} shows
instances of non-preserved, preserved, and strictly preserved faces. We will
generally drop the ``strictly'' whenever we talk about strictly preserved
faces.

What makes strictly preserved faces so nice is the fact that the above
conditions can be checked prior to the projection by purely (linear) algebraic
means. The key to that is the Projection Lemma (cf.  \cite[Proposition
3.2]{zie04}). We will give a variant of it which requires a few more notions. 

Let $P \subset \R^n$ be a full-dimensional polytope with $0$ in the interior.
The dual polytope is
\[
    P^\dual \ =\  \{ \ell \in \Rd{n} : \ell(x) \leq 1\,\, \textrm{ for all } x
    \in P \} = \conv \{ \ell_1, \dots, \ell_m \} \ \subset\  \Rd{n}
\]
and for every face $F\subseteq P$ we denote by $F^\fdual = \{ \ell \in P^\dual : \ell|_F = 1 \}$ the corresponding face 
of $P^\dual$. Furthermore, we define
$I : \{\text{faces of } P\} \rightarrow 2^{[m]}$ to be the map satisfying 
\[
     \conv \{ \ell_i : i \in I(F) \}  = F^\fdual 
\]
for every face $F \subseteq P$.  Let $\pi : \R^n \rightarrow \R^d$ be a
linear projection and let $q$ be a map fitting into the short exact sequence
\[
    0 \longrightarrow \R^{n-d} \stackrel{q}{\longrightarrow} \R^n
    \stackrel{\pi}{\longrightarrow} \R^d \longrightarrow 0.
\]
Dualizing gives rise to a (dual) exact sequence 
\[
    0 \longleftarrow \Rd{n-d} \stackrel{q^*}{\longleftarrow} \Rd{n}
    \stackrel{\pi^*}{\longleftarrow} \Rd{d} \longleftarrow 0.
\]
The characterization of strictly preserved faces will be in terms of the dual
map $q^*$ and the dual to the face under consideration. 

\begin{lem}[Projection Lemma]
    Let $P$ be a polytope and $F\subset P$ a face. Then $F$ is strictly
    preserved iff $0 \in \int q^*(F^\fdual) = \int \conv\{ q^*(\ell_i) : i \in
    I(F)\}$.
\end{lem}

We first sort out the situation for (non-strictly) preserved faces. 

\begin{prop}\label{prop:proj_is_face}
    Let $F \subset P$ be a face. Then $\pi(F)$ is a face of $Q = \pi(P)$ iff
    $0 \in q^*(F^\fdual)$.
\end{prop}

\begin{proof}
    $p(F)$ is a face of $Q$ iff there is an $\ell \in \Rd{d}$  such that
    \[
        \begin{array}{rl}
            \ell\circ \pi(x) < 1 & \text{for all } x \in P \backslash F\text{ 
            and} \\
            \ell\circ \pi(x) = 1 & \text{for all } x \in F \\
        \end{array}
    \]
    This is the case iff $\pi^*(\ell) \in F^\fdual$ which in turn holds iff $0 =
    q^* \circ \pi^* (\ell) \in q^*(F^\fdual)$.
\end{proof}

\begin{prop}\label{prop:proj_is_comb_equiv}
    Let $F\subset P$ be a face. Then $\pi(F)$ is combinatorially equivalent to
    $F$ iff $q^*(F^\fdual)$ is full-dimensional.
\end{prop}
\begin{proof}
    Let $\aff F = x_0 + L$ with $L$ being a linear space. Then $\aff F^\fdual =
    \ell_0 + L^0$ with $L^0 = \ker( \Rd{n} \rightarrow L^*)$.  The result
    follows if we show that $\pi$ is injective on $L$ iff $q^*$ restricted to
    $L^0$ is surjective. Now $\pi|_L$ is injective iff $\im q \cap L = 0$.
    This holds iff $q^*$ is surjective when restricted to $L^0 = (\R^n/L)^*$.
\end{proof}

\begin{proof}[Proof of Projection Lemma]
    Let $F \subset P$ be a preserved face, i.e.\ $\pi(F)$ is a face of $Q =
    \pi(P)$ combinatorially equivalent to $F$, and let $G = \pi^{-1}(\pi(F))$.
    Clearly, $G$ is a face of $P$ and $F \subseteq G$. 

    The inclusion is strict iff $G^\fdual$ is a proper face of $F^\fdual$.  Now,
    $\pi(G)$ is a face iff $0 \in q^*(G^\fdual)$ by Proposition
    \ref{prop:proj_is_face}. This, in turn, is the case iff $0 \not\in \int
    q^*(F^\fdual)$.
\end{proof}

\subsection*{The geometric side}

We made use of the fact that $q^*(F^\fdual) = \conv\{ q^*(w) : w \in \V
F^\fdual\}$ and for later reference we denote by $G = \{ g_i = q^*(\ell_i) : i
= 1,\dots,m\}$ the projection of the vertices of $P^\dual$.  Note that we will
not treat $G$ as the set of vertices of a polytope (especially since not all
would be vertices) but as a configuration of vectors. In case that all
vertices survive the projection, this vector configuration has some strong
properties: it is a \emph{Gale transform}. Gale transforms are a well-known
notion from discrete geometry; we refer the reader to Matou\v{s}ek
\cite{mat02} and Ziegler \cite{zie95} for full treatments (from different
perspectives) and McMullen \cite{mcm79} for an extensive survey.

A set of vectors $W = \{ w_1, \dots, w_k \} \subset \R^d$ is \emph{positively
spanning} if every point in $\R^d$ is a non-negative combination of the
vectors $w_i$, that is, if $\cone W = \R^d$. Equivalently, $U$ is positively
spanning if $\conv W$ is a full dimensional polytope with $0$ in its interior.
We also need the weaker notion of \emph{positively dependent} which holds if
$0 \in \mathsf{relint}\, \conv W$.

\begin{dfn}[Gale transform] 
    A finite vector configuration $G = \{g_1,\dots,g_m\}$ 
    is a \emph{Gale transform} if for every $i = 1,\dots,m$ the
    subconfiguration $G\backslash g_i$ is positively spanning.
\end{dfn}

The main reason why Gale transforms are useful is that they are yet another
way to represent polytopes.
\begin{thm}[Gale duality]\label{thm:gale_duality}
    Let $G = \{g_1,\dots,g_m\}$ be a Gale transform in $(m-d-1)$-dimensional
    space, then there is a $d$-polytope $Q$ with vertices $V =
    \{v_1,\dots,v_m\}$ such that for every $I \subset [m]$ 
    \[
        \conv\{v_i : i \in [m]\backslash I\} \text{ is a face of $Q$} \ \Longleftrightarrow\ 
        \{g_j : j \in I \} \text{ are positively dependent}.
    \]
    Furthermore, $Q$ is unique up to affine isomorphisms.
\end{thm}

The condition given by the Projection Lemma can be rephrased in terms of
positive spans. The key observation is that the set $G$ is actually a Gale
transform if all vertices survive the projection.

\begin{prop}\label{prop:is_gale}
    Let $P \subset \R^n$ be a $n$-polytope with $m$ facets and $\pi: \R^n
    \rightarrow \R^d$ a linear projection.  The set $G$ is a Gale transform if
    all vertices of $P$ are strictly preserved under $\pi$.
\end{prop}
\begin{proof}
    Since $P$ is full dimensional for every $i \in [m]$ there is a vertex $v
    \in \V P$ such that $i \not\in I(v)$. Since $v$ is preserved under
    projection the set $\{g_j : j \in I(v) \} \subset G\backslash g_i$ is
    positively spanning. 
\end{proof}

Thus, by Theorem \ref{thm:gale_duality}, there is a combinatorially unique
polytope $\A(P,\pi)$ associated to $(P,\pi)$ which we simply call the
\emph{associated polytope}. So far, the main property of $\A(P,\pi)$ is that
it sort of witnesses the survival of the vertices.

For the case that interests us, namely $P$ being a product of simplices, we
can even assume that $\A(P,\pi)$ is a simplicial polytope, the reason being
the following: If $P^\Delta$ is simplicial we can wiggle the vertices without
changing the combinatorial type. For strictly preserved faces of $P$, the
condition as dictated by the Projection Lemma is \emph{open}, i.e.\ stable
under small perturbations. Thus perturbing the bounding hyperplanes of $P$
does neither alter the type nor the fact that all vertices survive the
projection.  The effect on the Gale transform $G$ is that the vectors are in
general position with respect to hyperplanes containing the origin. On the
other hand, this is yet another characterization of the fact that $\A(P,\pi)$
is a polytope with vertices in general positions and hence simplicial.

Summing up so far, we have the following.
\begin{cor} \label{cor:associated_to_simplex_product}
    Let $P \subset \R^{d^2}$ be a polytope combinatorially equivalent to a
    $d$-fold product of $d$-simplices such that a projection to $d$-space
    preserves all the vertices. Then there is a $(2d-1)$-dimensional,
    simplicial polytope $\A(P,\pi)$ with $d(d+1)$ vertices associated to the
    projection.
\end{cor}

\subsection*{The combinatorial side}

Given that all vertices of a polytope $P$ survive the projection we obtain an
associated polytope $\A = \A(P,\pi)$. Furthermore, for every vertex $v \in P$
the polytope $q^*(v^\fdual)$ has zero in its interior which, in particular,
implies that its vertices  are positively dependent. By Gale duality (Theorem
\ref{thm:gale_duality}), this induces a face in $\A$. The collection of all
the induced faces is a polytopal complex in the boundary of $\A$. If $P$ is a
simple polytope, we argued that $\A$ can be assumed to be simplicial. Thus,
the polytopal complex is a simplicial complex whose combinatorics is
determined by the sole knowledge of the combinatorics of $P$. We give a rather
general description of this complex since it seems that it is the first
occurrence in the literature.  For background on simplicial complexes as well as
notation, we refer to the first chapter of \cite{MatousekBZ:BU}.

\begin{dfn}[Complement complex] 
    Let $\K \subseteq 2^V$ be a simplicial complex on vertices $V$. The
    \emph{complement complex} $\cc{\K}$ of $\K$ is the closure of
    \[
        \{ V \backslash \tau : \tau \in \K, \ \tau \text{ a facet}\}.
    \]
\end{dfn}

From the bare definition of the complement complex we deduce the following
simple properties whose proofs we omit.

\begin{prop}
    Let $\K$ and $\Ll$ be simplicial complexes. Then the following statements
    hold
    \begin{enumerate}
        \item[\rm 1.] $\cc{(\cc{\K})} = \K$
        \item[\rm 2.] $\dim \cc{\K} = n - \dim\K - 1$.
        \item[\rm 3.] $\cc{(\K * \Ll)} = \cc{\K} * \cc{\Ll}$.
    \end{enumerate}
\end{prop}

In particular, the first property states that no information is lost in the
passage from $\K$ to its complement complex.  To the best of our knowledge
there has been no work on this construction of simplicial complexes. A
possible reason for that is the lack of topological plausibility. For a
simplicial complex $\K$ there seem to be no obvious relation between $\K$ and
$\cc{\K}$ concerning homotopy type and/or (co)homology.  For a complex being a
\emph{matroid}, the complement complex is just the matroid dual which is well
understood combinatorially and topologically. But matroids are rare among
simplicial complexes. 

For a simplicial polytope $P^\dual$ we denote by $\bd(P^\dual)$ the simplicial
complex of all proper faces of $P^\dual$.  For a fixed numbering of the facets
of $P$, the vertex set of $\bd(P^\dual)$ can be identified with $[m] = \{1,
\dots, m \}$.

\begin{thm}
    Let $P$ be a simple polytope whose vertices are preserved under $\pi$.
    Then the complex $\cc{\bd(P^\dual)}$ is realized in the
    boundary of the associated polytope $\A(P,\pi)$.
\end{thm}

\begin{proof}
    Since $P$ is simple, the associated polytope $\A(P,\pi)$ is simplicial.
    Therefore, it is sufficient to show that all facets of $\cc{\bd(P^\dual)}$
    are part of the boundary. For now, let $q_1,\dots,q_m$ be the vertices of
    $\A(P,\pi)$ labelled in accordance with the elements of its Gale
    transform.

    A facet of $\cc{\bd(P^\dual)}$ is of the form
    $[m] \backslash I(v)$ for some vertex $v \in \V P$. Since
    $\pi(v)$ is a vertex of the projection, we have $0 \in \int \conv \{ g_i :
    i \in I(v)\}$ by the Projection Lemma. By Gale duality, this corresponds
    to the fact that $\conv\{q_i : i \in [m]\backslash I(v)\}$ is a face of 
    $\A(P,\pi)$.
\end{proof}

For a full-dimensional polytope $P \times Q$ with $0$ in its interior the
polar dual is $P^\dual \oplus Q^\dual$ whose proper faces are the joins of
proper faces of $P^\dual$ and of $Q^\dual$. Thus for the complement of the
boundary complex we have
\[
    \cc{\bd\left( (P \times Q)^\dual \right)} = \cc{\bd(P^\dual)} *
    \cc{\bd(Q^\dual)}.
\]

Focusing again on the polytopes in question, we want to consider the
complement complex for the dual of a $d$-fold product of simplices. Since a
simplex is self-dual, we have that $\bd(\Delta_d) \cong \tbinom{ [d+1]}{\le
d}$. Thus, for a $d$-fold product of $d$-simplices the corresponding
complement complex is equivalent to 
\[
    \tbinom{[d+1]}{ \le 1}^{*d},
\]
that is, the $d$-fold join of a complex consisting of $d+1$ isolated points.

\begin{cor}\label{cor:bla} 
    If there is a realization of a $d$-fold product of $d$-simplices such that
    a projection to $d$-space retains all vertices, then the complex
    $\tbinom{[d+1]}{\le 1}^{*d}$ is embeddable into a sphere of dimension
    $2d-2$.
\end{cor}

This will be our punchline: We will show that the embedding claimed by
Corollary \ref{cor:bla} does not exist. Let us rest for a moment and
reconsider the example from Section \ref{sec:problem_reductions}.

Let $P$ be a realization of a product of two triangles such that a projection
$\pi$ to the plane preserves all (nine) vertices. By Corollary
\ref{cor:associated_to_simplex_product}, the associated polytope $\A(P,\pi)$
is a $3$-dimensional simplicial polytope with $6$ vertices. The complement
complex is $\tbinom{[3]}{\le 1}^{*2}$, can be thought of as $\K_{3,3}$, the
complete bipartite graph on 6 vertices. By Corollary \ref{cor:bla}, this
complex is embedded in the boundary of $\A(P,\pi)$, which is a $2$-sphere.
This, however, is impossible: Graphs embeddable into the $2$-sphere are planar
while the $\K_{3,3}$ is minimal non-planar.

For a $3$-fold product of $3$-simplices the boundary of the associated
polytope is a $4$-sphere and the complement complex is $2$-dimensional. So,
again, there is no (obvious) elementary reasoning neither are there
off-the-shelf results showing non-embeddability. We therefore have to resort
to more sophisticated machinery, as presented in the following section.

\section{Interlude: Embeddability of simplicial complexes}
\label{sec:embeddability}

We only give an \emph{executive summary} of the techniques and results
needed for the following; see \cite{MatousekBZ:BU}. 

The category of free $\Zt$-spaces consists of topological spaces $X$ together
with a free action of the group $\Zt$, i.e.\ a fixed point free involution on
$X$.  Morphisms in this category are continuous maps that commute with the
respective $\Zt$-actions. The foremost examples of $\Zt$-spaces are spheres
$\S^d$ with the antipodal action. For a $\Zt$-space $X$ a numerical invariant
is the \emph{$\Zt$-index} $\ind X$ which is the smallest integer $d$ such that
there is a $\Zt$-equivariant map $X \Ztmap \S^d$. For example  $\ind \S^d =
d$, which is a statement equivalent to the Borsuk--Ulam theorem. 

For a simplicial complex $\K$ we define the \emph{deleted join} of $\K$ to be
the complex
\[
    \djn{\K} = \{ \sigma \uplus \tau : \sigma, \tau \in \K, \sigma \cap \tau =
    \emptyset \}
\]
The deleted join turns an arbitrary simplicial complex into a free
$\Zt$-complex by means of $\sigma \uplus \tau \mapsto \tau \uplus \sigma$.

\begin{thm}[\cite{MatousekBZ:BU}, Theorem 5.5.5]\label{thm:Z2_embedd}
    Let $\K$ be a simplicial complex. If
    \[
        \ind \djn{\K} > d,
    \]
    then $\K$ is not embeddable into the $d$-sphere.
\end{thm}

The $\Zt$-index is rather difficult to calculate for general spaces. Luckily,
for the situation in which we will apply Theorem \ref{thm:Z2_embedd} there is
a beautiful theorem due to Karanbir Sarkaria (see \cite{MatousekBZ:BU}). In
order to state it properly we need some more definitions.

{\bf Minimal non-faces.} Let $\K \subset 2^V$ be a simplicial complex. A set
$F \subset V$ is called a \emph{non-face} if $F \not\in \K$ and its a
\emph{minimal non-face} if every proper subset of $F$ is in $\K$. We will
denote by $\nf(\K)$ the set of minimal non-faces.

{\bf Generalized Kneser graphs.} For a collection of sets $\nf = \{ F_1,
\dots, F_k\}$ we denote by $\KG(\nf)$ the (abstract) graph with vertex set
$\nf$. Two vertices $F_i, F_j$ share an edge iff $F_i \cap F_j = \emptyset$.
Such a graph is called a \emph{generalized Kneser graph}.

Finally, for a graph $G$ we denote by $\chi(G)$ its \emph{chromatic number},
i.e.\ the minimal number of colors to properly color the graph.

\begin{thm}[Sarkaria's coloring/embedding theorem]\label{thm:sarkaria}
    Let $\K$ be a simplicial complex with $n$ vertices and let $\nf = \nf(\K)$
    be the set of minimal non-faces. Then
    \[
        \ind \djn{\K} \ge n - \chi(\KG(\nf)) - 1.
    \]
\end{thm}

Taking up the example of triangle times triangle for the last time, let us
use Theorem \ref{thm:sarkaria} to show that $\K = \tbinom{[3]}{\le 1 }^{*2}$
does not embed into the $2$-sphere. This complex has $6$ vertices and, for
reasons we will give in the next section, the Kneser graph of its non-faces
is, again, the complete bipartite graph $\K_{3,3}$. Thus, using Sarkaria's
theorem, we get
\[
    \ind \djn{\K} \ge 6 - 2 - 1 = 3,
\]
which shows that $\K_{3,3}$ is not planar.

\section{Analysis of the complement complex}
\label{sec:compl_complex}

Although determining upper bounds on the chromatic number of graphs is easier
than finding equivariant maps, it is, in general, still hard enough. The key
property that enables us to calculate chromatic numbers for the Kneser graphs
we will encounter is that the complexes are made up of (possibly) simpler
ones, that is they are joins of complexes. The following results will show
that this continues to hold if we pass from complexes to non-faces and then to
Kneser graphs.

\begin{lem}\label{lem:nf_of_join}
    Let $\K$ and $\Ll$ be simplicial complexes. Then 
    \[
        \nf(\K * \Ll) = \{ F \uplus \emptyset : F \in \nf(\K) \}  \cup
                   \{ \emptyset \uplus G : G \in \nf(\Ll) \}.
    \]
\end{lem}
\begin{proof}
    Let $F \uplus G \in \nf(\K * \Ll)$ and $i \in F$ and $j \in G$.  Since $F
    \uplus G$ is a minimal non-face, it follows that $F\backslash i \uplus G$
    and $F \uplus G\backslash j$ are both in $\K*\Ll$. This, however, implies
    that $F \in \K$ and $G \in \Ll$.
\end{proof}

On the level of Kneser graphs this fact results in a \emph{bipartite sum} of
the respective Kneser graphs. Let $G$ and $H$ be graphs with disjoint vertex
sets $U$ and $V$. The bipartite sum of $G$ and $H$ is the graph $G \bowtie H$
with vertex set $U \cup V$ and edges $E(G) \cup E(H) \cup (U \times V)$. 

\begin{prop} \label{prop:chi_sum}
    Let $G$ and $H$ be graphs. Then
    \[
        \chi(G \bowtie H) = \chi(G) + \chi(H).
    \]
\end{prop}
\begin{proof}
    The edges $U \times V \subset E(G \bowtie H)$ force the set of colors on
    $U$ and $V$ to be disjoint. Thus a coloring on $G \bowtie H$ is minimal
    iff it is minimal on the subgraphs $G$ and $H$.
\end{proof}

The complex to which we want to apply the result is $\K = \Ll^{*d}$ with $\Ll =
\tbinom{[d+1]}{\le 1}$. We will analyze the chromatic number of $\KG(\nf)$ for
$\nf = \nf(\Ll)$ and use Proposition \ref{prop:chi_sum} to get an obstruction
to the embeddability of $\K$ into some sphere.

From the definition of $\Ll$ we see that the minimal non-faces are exactly the
two element subsets of $[d+1]$, that is $\nf = \tbinom{[d+1]}{2}$. The
resulting Kneser graph $\KG(\nf)$ is an instance of a famous family of graphs,
the {\it ordinary} Kneser graphs $\KG_{n,k} := \KG\tbinom{[n]}{k}$. The
determination of their chromatic numbers is one of the first success stories
of topological combinatorics.

\begin{thm}[Lov\'{a}sz--Kneser theorem \cite{lov78}] 
    For $0 < 2k - 1 \le n$ the chromatic number of the Kneser graph
    $\KG_{n,k}$ is $\chi(\KG_{n,k}) = n - 2k + 2$.
\end{thm}

With that last bit of information we can finally complete the proof of Theorem
\ref{thm:main_stronger}.

\begin{proof}[Proof of Theorem \ref{thm:main_stronger}]
    Assume that there is a realization of a $d$-fold product of $d$-simplices
    whose projection to $d$-space preserves all the vertices. This implies, by
    Corollary \ref{cor:bla}, that the complex $\K = \tbinom{[d+1]}{\le 1
    }^{*d}$ is embeddable into a $(2d-2)$-sphere.

    Let $\Ll = \tbinom{[d+1]}{\le 1}$ and let $\nf(\Ll) = \tbinom{[d+1]}{2}$ be
    the minimal non-faces of $\Ll$. For the associated Kneser graph
    $\KG(\nf(\Ll)) = \KG_{d+1,2}$ we have $\chi(\KG_{d+1,2}) = d-1$.

    By Lemma \ref{lem:nf_of_join} and Proposition \ref{prop:chi_sum} we have for
    $\nf = \nf(\K)$
    \[
        \ind \djn{\K} \ge d(d+1) - \chi( \KG(\nf)) - 1 = 
        d(d+1) - d \chi(\KG_{d+1},2) - 1 =  2d -1 > 2d -2
    \]
    By Theorem \ref{thm:Z2_embedd}, this contradicts the claim that $\K$ is
    embeddable into a $(2d-2)$-sphere.
\end{proof}

It is true that any upper bound on the chromatic number of $\KG_{d+1,2}$ would
have sufficed and it was thus unnecessary to invoke the Lov\'{a}sz--Kneser
theorem. We, nevertheless, wish to argue that the application of the
Lov\'{a}sz--Kneser is justified. Sarkaria's theorem can be used with any upper
bound on the chromatic number of the Kneser graph. This, however, results in a
weaker bound on the $\Zt$-index of $\djn{\K}$. Using the actual chromatic
number shows that Theorem \ref{thm:main_stronger} is sharp concerning the
number of factors, i.e. there are no topological obstruction for a product of
less than $d$ simplices. On the other hand, by Proposition 5.3.2 in \cite[p.
96]{MatousekBZ:BU}, we have that $2d-1 \ge \ind \djn{\K}$ and, therefore, the
calculation in the preceding proof gives the actual $\Zt$-index. Thus, Theorem
\ref{thm:main_stronger} is also sharp with respect to the dimension of the
target space, i.e.\ projecting to a space of dimension $\ge d$. This, in
particular, stands in favor for the result of Fukuda \& Weibel.

\section{Remarks}\label{sec:remarks}

At the Oberwolfach-Workshop ``Geometric and Topological Combinatorics'' in
January 2007 Rade \v{Z}ivaljevi\'{c} suggested a different argument involving
Lov\'{a}sz' \emph{colored Helly theorem}.

\begin{thm}[Colored Helly Theorem; cf.\ \cite{lov74}] 
    Let $\C_1,\dots,\C_r$ be collections of convex sets in $\R^d$ with $r \ge
    d+1$.  If $\cap_{i=1}^r C_i \not= \emptyset$  for every choice $C_i \in
    \C_i$ then there is a $j \in [r] $ such that $\cap_{C \in \C_j} C \not=
    \emptyset$.
\end{thm}

The following proof was supplied by Imre B\'{a}r\'{a}ny (personal
communication): Let $P_1,\dots,P_r \subset \R^d$ be $d$-polytopes and let
$\C_i = \{C_v \subset \Rd{d} : v \in \V P_i\}$. The $C_v$ are defined by the
condition that $\ell \in C_v$ if and only if $\ell$ attains its unique maximum
over $P_i$ in $v$.  It is clear that the $C_v$ are pairwise disjoint. Now, if
$P_1 + \cdots + P_r$ attains the trivial upper bound then for every choice of
vertices $v_i \in \V P_i$ the intersection $C_{v_1} \cap \cdots \cap C_{v_r}$
is non-empty and, thus, contradicts the colored Helly theorem. 
In the colored Helly theorem the bound $r \ge d+1$ is tight and, thus, yields
Theorem \ref{thm:main} in a slightly weaker version with at least $d+1$
summands. 

\begin{figure}[ht]
    \centering
    \includegraphics[height=8cm]{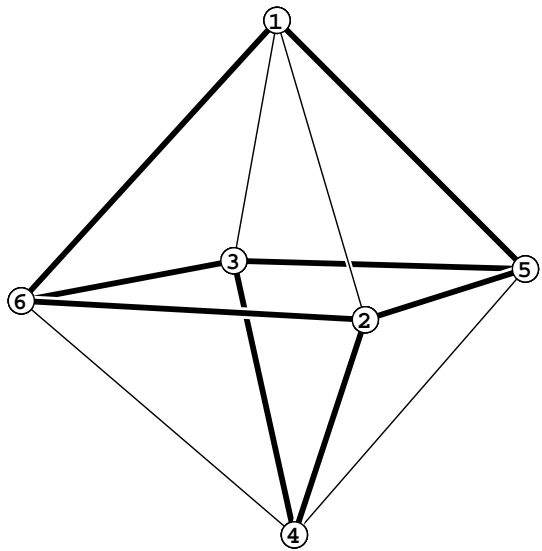}\\
    \caption{The polytope $\A$ associated to $P$. The marked edges correspond
    to the embedding of $\K - \{1,4\}$, that is
    $\K_{3,3}$ minus an edge.} \label{fig:octahedron}
\end{figure}

In Section \ref{sec:problem_reductions} we claimed the existence of a
combinatorial product of two triangles that projects to a plane $8$-gon.
Consider the polytope $P \subset \R^4$ given as the set of solutions to the
following system of inequalities
\[
    \begin{array}{r} 
    \scriptstyle   1 \\
    \scriptstyle   2 \\
    \scriptstyle   3 \\
    \scriptstyle   4 \\
    \scriptstyle   5 \\
    \scriptstyle   6 \\
    \end{array}
    \left(
    \begin{array}{rrrr} 
   1 &         1 &           &    \\ 
  -1 &         1 &           &    \\ 
   0 &        -1 &     -\eps &    \\ 
     &     -\eps &        -1 &  0 \\ 
     &           &         1 &  1 \\ 
     &           &         1 & -1 \\ 
        \end{array}\right) \mathbf{x} \le 
    \left( \begin{array}{r} 
       1 \\ 1 \\ 1 \\ 1 \\ 1 \\ 1 \\
   \end{array}\right)
\]        
The numbers to the left label the facets of $P$.  For $\eps = 0$ this is just
a Cartesian product of two triangles and, since this is a simple polytope, we
can choose $\eps > 0$ sufficiently small without changing the combinatorial
type. Taking $\pi : \R^4 \rightarrow \R^2$ to be the projection to the first
and last coordinate and identifying $\Rd{4}$ with  $\R^4$ via the standard
scalar product we get the ordered set
\[
    G = q^*(\V P^\dual) = \left(
    \begin{array}{rrrrrr} 
          1 &  1 &    -1 & -\eps &    &    \\
            &    & -\eps &    -1 &  1 &  1 \\
    \end{array}\right)
\]        
For $0 < \eps \le 1$ the set $G$ is the Gale transform of a polytope $\A$
combinatorially equivalent to an octahedron (e.g.\ set $\eps = 1$ and observe
that $G$ is a Gale transform of a regular octahedron.) As intersections of
facets the set of vertices is given by $S = \{ [6] \backslash \{i,j\} :
\{i,j\} \in \K  \}$ where the complement complex $\K$ is the complete
bipartite graph on the partition $\{1,2,3\}$ and $\{4,5,6\}$.
We show that the only vertex $v_0$ that fails to
survive the projection is given by the intersection of the facets $[6]
\backslash \{1,4\}$.  By the Projection Lemma and Gale duality this is case if
and only if $\K - \{1,4\}$ is a subcomplex of the $1$-skeleton of
$\A$. Figure \ref{fig:octahedron} shows $\A$ and the embedding of $\K -
\{1,4\}$, thus finishing the proof.
The missing edge between the vertices $1$ and $4$ shows that $v_0$ fails
short of being a vertex of $\pi(P)$.

\bibliographystyle{siam}
\begin{small}
    \bibliography{/Users/sanyal/Documents/Bibliography,AddtionalBib}

\begin{thebibliography}{10}

\bibitem{am98}
{\sc N.~Amenta and G.~M. Ziegler}, {\em Deformed products and maximal shadows
  of polytopes}, in Advances in Discrete and Computational Geometry, R.~P.
  B.~Chazelle, J.E.~Goodman, ed., vol.~223 of Contemporary Mathematics, Amer.
  Math. Soc., Providence RI, 1998, pp.~57--90.

\bibitem{bar07}
{\sc I.~B\'{a}r\'{a}ny}, {\em email communication}.
\newblock Janurary 2007.

\bibitem{DBLP:journals/cad/ElberK99}
{\sc G.~Elber and M.-S. Kim}, {\em Offsets, sweeps, and {M}inkowski sums},
  Computer-Aided Design, 31 (1999).

\bibitem{fuk06}
{\sc K.~Fukuda and C.~Weibel}, {\em On $f$-vectors of {M}inkowski additions of
  convex polytopes}, Discrete \& Computational Geometry, to appear.

\bibitem{Goldfarb2}
{\sc D.~Goldfarb}, {\em On the complexity of the {simplex} method}, in Advances
  in optimization and numerical analysis (Oaxaca, 1992), S.~Gomez, ed.,
  Dordrecht, 1994, Kluwer, pp.~25--38.
\newblock Proc.\ 6th Workshop on Optimization and Numerical Analysis, Oaxaca,
  Mexico, January 1992.

\bibitem{latombe}
{\sc J.-C. Latombe}, {\em Robot {M}otion Planning}, vol.~124 of Series in
  Engineering and Computer Science, Kluwer International, 1991.

\bibitem{lov74}
{\sc L.~Lov\'{a}sz}, {\em Problem 206}, Matematikai Lapok, 25 (1974), p.~181.

\bibitem{lov78}
\leavevmode\vrule height 2pt depth -1.6pt width 23pt, {\em Kneser's conjecture,
  chromatic number, and homotopy}, J. Combinatorial Theory, Ser.~A, 25 (1978),
  pp.~319--324.

\bibitem{mat02}
{\sc J.~Matou{\v s}ek}, {\em Lectures on Discrete Geometry}, vol.~212 of
  Graduate Texts in Mathematics, Springer-Verlag, New York, 2002.

\bibitem{MatousekBZ:BU}
\leavevmode\vrule height 2pt depth -1.6pt width 23pt, {\em {Using the
  {Borsuk--Ulam} Theorem. {L}ectures on Topological Methods in Combinatorics
  and Geometry}}, Universitext, Springer-Verlag, Heidelberg, 2003.
\newblock Written in cooperation with Anders Bj\"orner and G\"unter M. Ziegler.

\bibitem{mcm79}
{\sc P.~McMullen}, {\em Transforms, diagrams and representations}, in
  Contributions to Geometry, J.~T{\"o}lke and J.~Wills, eds., Proc. Geom.
  Sympos., Siegen, Birkh{\"a}user Basel, 1979, pp.~92--130.

\bibitem{rsz07}
{\sc T.~R\"{o}rig, R.~Sanyal, and G.~M. Ziegler}, {\em Polytopes and surfaces
  via projection}, in preparation.

\bibitem{stu02}
{\sc B.~Sturmfels}, {\em Solving Systems of Polynomial Equations}, no.~97 in
  CBMS Regional Conferences Series, Amer. Math. Soc., Providence RI, 2002.

\bibitem{zie95}
{\sc G.~M. Ziegler}, {\em Lectures on Polytopes}, vol.~152 of Graduate Texts in
  Mathematics, Springer-Verlag, New York, 1995.

\bibitem{zie04}
\leavevmode\vrule height 2pt depth -1.6pt width 23pt, {\em Projected products
  of polygons}, Electron. Res. Announc. Amer. Math. Soc., 10 (2004),
  pp.~122--134.
\newblock \url{http://www.ams.org/era/2004-10-14/S1079-6762-04-00137-4}.

\end{thebibliography}
\end{small}

\end{document}